# Truth and Knowledge: the Incorrect Definition of 'Powers'
## A New Interpretation of *Theaetetus* (147d7-148b2)

Luc Brisson
Centre Jean Pépin
CNRS-UMR 8230

Salomon Ofman
Institut mathématique de Jussieu-
Paris Rive Gauche/HSM
Sorbonne Université et Université Paris 7

**Abstract.** In a first article (referred here as B-O), we studied the first part of the so-called 'mathematical part' of Plato's *Theaetetus*, i.e. Theodorus' lesson. In the present one, we consider the sequel and the end of the passage (147d7-148b2), as well as its philosophical interpretation in connection with the whole dialogue. As in the previous article, we analyze it simultaneously from the mathematical, the historical and the philosophical points of view, a necessity to understand it. Our strategy is once again to take seriously Plato's text, not as the dream of a poet.

Our analysis casts a new light on this passage, as an essential testimony for both Plato's philosophy and for history of mathematics. Plato's relation to Theodorus and Theaetetus is more complex than usually claimed; the search for a definition of knowledge conducted in a large part of the dialogue and even in some other dialogues is rooted in the mathematical passage; it needs a reevaluation of its connection to Euclid's *Elements*, in particular the proposition X.9, as well as a supposed Euclid's 'catastrophic' mistake (to quote Jean Itard) on incommensurability. In a nutshell, the passage has to be understood differently than in the usual interpretations gathered together under the collective so-called name 'Modern Standard Interpretation' (or MSI).

Both articles form a whole. They are both aimed to an audience without any particular mathematical background, and require only elementary mathematical knowledge, essentially of high school-level. Some more complex points are developed in an **Appendix** at the end of the article.

**Résumé.** Dans un premier article (noté ici par B-O), nous avons présenté une nouvelle traduction et une nouvelle interprétation de la leçon mathématique de Théodore sur les grandeurs irrationnelles. Dans celui-ci, nous analysons la suite et la fin du passage (147d7-148b2), ce qui nous conduit à donner une interprétation philosophique globale le situant dans le cadre de la totalité du dialogue. Comme dans le précédent article, nous l'étudions simultanément des points de vue mathématique, historique et philosophique, une nécessité pour le comprendre. Notre ligne directrice est là encore de prendre le texte platonicien au sérieux, non comme le rêve d'un poète.

Notre analyse apporte un éclairage nouveau sur ce passage, en tant que témoignage essentiel à la fois pour la philosophie de Platon et l'histoire des mathématiques. La relation de Platon à Théodore et à Théétète est plus complexe que ce que l'on soutient communément; dans le passage mathématique s'enracinent toutes les questions discutées non seulement dans le dialogue, mais aussi dans certains autres ; il est enfin nécessaire de réévaluer son rapport aux *Éléments* d'Euclide, en particulier à la proposition X.9, ainsi que la supposée erreur

'catastrophique' d'Euclide (pour citer Jean Itard) sur l'incommensurabilité. En bref, ce passage doit être compris différemment des interprétations usuelles regroupées sous le nom collectif de 'Interprétation moderne standard' (ou son acronyme anglais MSI).

Ces deux articles forment un tout. Ils se destinent tous deux à un public sans formation mathématique particulière, sa compréhension ne supposant que des connaissances mathématiques très élémentaires, essentiellement celles des premières années de collège. Certains points plus délicats sont développés dans une **Annexe** à la fin de l'article.

### I. *Theaetetus* 147d-148b

In this second part, Theaetetus gives an account of some work that he did with his young friend *Socrates* [1] following Theodorus' lesson on 'powers' [2] that they had attended.

#### 1. The text

This translation follows strictly the Greek text. The most important difficulties, noted in bold characters in the Greek text, are developed and explained in the notes following it.

| | |
|---|---|
| THEAETETUS<br>So the idea occurred to us that, since these powers were turning out to be quantitatively unlimited, we might try to put them together into one single thing, by which we could name all these powers. | ΘΕΑΙ.<br>ἡμῖν οὖν εἰσῆλθέ τι τοιοῦτον, ἐπειδὴ ἄπειροι τὸ πλῆθος αἱ δυνάμεις ἐφαίνοντο, πειραθῆναι συλλαβεῖν **εἰς ἕν**, ὅτῳ πάσας ταύτας προσαγορεύσομεν τὰς δυνάμεις. |
| SOCRATES<br>[**147e**] And then did you find such a thing? | ΣΩ.<br>Ἦ καὶ ηὕρετέ τι τοιοῦτον; |
| THEAETETUS<br>I think we did. But consider if you agree too. | ΘΕΑΙ.<br>Ἔμοιγε δοκοῦμεν· σκόπει δὲ καὶ σύ. |
| SOCRATES<br>Go on. | ΣΩ.<br>Λέγε. |
| THEAETETUS<br>We divided the integer in its totality in two parts. One which has the power to be a product of an equal times an equal [**148a**], we likened it to the square as a figure, and we named it square and equilateral. | ΘΕΑΙ.<br>**Τὸν ἀριθμὸν πάντα** δίχα διελάβομεν· τὸν μὲν **δυνάμενον ἴσον ἰσάκις** γίγνεσθαι τῷ **τετραγώνῳ** τὸ σχῆμα ἀπεικάσαντες **τετράγωνόν** τε καὶ ἰσόπλευρον προσείπομεν. |
| SOCRATES<br>Good, so far. | ΣΩ.<br>Καὶ εὖ γε. |
| THEAETETUS<br>Then, the one in-between, such as three | ΘΕΑΙ.<br>Τὸν τοίνυν μεταξὺ τούτου, ὧν καὶ τὰ |
| or five or any one which has not the power to be a product of an equal times an equal, but is the product of a greater times a lesser, or a lesser times a greater, that is, always encompassed by a greater side and a smaller one, we likened it this time to a rectangular figure so that we called it a rectangular integer. | τρία καὶ τὰ πέντε καὶ πᾶς ὃς ἀδύνατος ἴσος ἰσάκις γενέσθαι, **ἀλλ' ἢ πλείων ἐλαττονάκις ἢ ἐλάττων πλεονάκις γίγνεται, μείζων δὲ καὶ ἐλάττων ἀεὶ πλευρὰ αὐτὸν περιλαμβάνει**, τῷ προμήκει αὖ σχήματι ἀπεικάσαντες προμήκη ἀριθμὸν ἐκαλέσαμεν. |
| SOCRATES<br>That's excellent. But how did you go on?<br>THEAETETUS<br>All the lines squaring the equilateral and plane integer, we defined as 'length', while we defined as 'powers' the ones squaring the rectangular, because, although not commensurable as length with the formers, [**148b**] they are commensurable as areas that they have the power to produce. And it is the same in the case of solids.<br>SOCRATES<br>Excellent my boys. I don't think Theodorus is likely to be up for false witness. | ΣΩ.<br>Κάλλιστα. ἀλλὰ τί τὸ μετὰ τοῦτο;<br>ΘΕΑΙ.<br>Ὅσαι μὲν **γραμμαὶ** τὸν ἰσόπλευρον καὶ ἐπίπεδον ἀριθμὸν τετραγωνίζουσι, μῆκος ὡρισάμεθα, ὅσαι δὲ **τὸν ἑτερομήκη**, δυνάμεις, ὡς μήκει μὲν **οὐ συμμέτρους ἐκείναις, τοῖς δ'ἐπιπέδοις ἃ δύνανται. Καὶ περὶ τὰ στερεὰ ἄλλο τοιοῦτον**.<br><br>ΣΩ.<br>**Ἄριστά γ' ἀνθρώπων**, ὦ παῖδες· ὥστε μοι δοκεῖ ὁ Θεόδωρος οὐκ ἔνοχος τοῖς ψευδομαρτυρίοις ἔσεσθαι. |

**2. Commentary**

The text is difficult to understand because history, philosophy, philology and mathematics are closely intertwined here. Moreover, there is almost no extant text about Greek mathematics from the end of the 5th century, let alone from an earlier period. As a matter of fact, the mathematical passage is one of the most important textual testimonies about this period. Hence, we have to reconstruct more or less from scratch the mathematics underlain by Theaetetus' account and simultaneously to avoid the dangerous path of anachronisms. [3]

1) 'one single thing' ('εἰς ἕν'). On trouve 49 fois cette expression dans le corpus platonicien. Les emplois les plus proches se trouvent dans le *Sophiste*. Le premier porte sur l'ensemble que constituent les simulacres produit pas le sophiste que l'Étranger d'Élée vient de ranger sous la technique de l'apparence. Théétète répond : 'Aucune. Tu viens de parler d'une espèce très nombreuse ('πάμπολυ γὰρ εἴρηκας εἶδος'), en rassemblant tout dans une seule chose ('εἶδος εἰς ἓν πάντα συλλαβὼν') ; c'est, il me semble, pratiquement la forme la plus variée' (234b3-4). Plus loin l'Étranger d'Élée explique : 'En outre, ceux qui tantôt unissent toutes choses et tantôt les divisent, soit qu'ils les amènent à l'unité ou que, de l'unité ils fassent sortir une infinité de choses ('εἰς ἓν καὶ ἐξ ἑνὸς ἄπειρα εἴτε εἰς πέρας ἔχοντα στοιχεῖα διαιρούμενοι καὶ ἐκ τούτων συντιθέντες'), soit qu'ils les divisent jusqu'à des éléments qui ont une limite et, à partir de ceux-ci, produisent une unité, pareillement lorsqu'ils admettent que cela arrive périodiquement ou éternellement, dans tous les cas, ils ne diraient rien si toutefois il n'y a aucun mélange.' (252b1-6) Plus loin encore, l'Étranger d'Élée fait cette déclaration qui se rapproche de notre texte : 'Conservons donc cette partie sous le nom de 'genre mimétique'. Laissons tomber, en revanche, l'autre partie, non seulement par paresse, mais aussi afin de permettre que quelqu'un d'autre trouve un jour son unité et lui applique un nom approprié ('τὸ δ' ἄλλο πᾶν ἀφῶμεν μαλακισθέντες καὶ παρέντες ἑτέρῳ συναγαγεῖν τε εἰς ἓν καὶ πρέπουσαν ἐπωνυμίαν ἀποδοῦναί τιν' αὐτῷ.'). (267a10-b2). L'idée est toujours la même. Il s'agit de tout réunir sous un seul référent auquel on donne le nom 'sophiste' objet de la recherche dans le dialogue du même nom. Et il en va de même dans le *Politique* pour le 'politique' (267b5-c1; 308c1-7).

2) 'the integer in its totality' ('**Τὸν ἀριθμὸν πάντα**'). Il y a ambiguïté. On peut en effet penser aux entiers en leur totalité (c'est la totalité qui est divisé en deux) ou à chaque entier en particulier (c'est chacun d'entre eux qui est divisé en deux). Il n'est aucun doute, d'après la suite, qu'il s'agit ici du premier cas. Voir également *Phèdre* 245c5 pour une telle ambiguïté.

3) 'the power to be an equal times an equal' ('**δυνάμενον ἴσον ἰσάκις**'). The Greek text says 'equally (ἰσάκις) equal (ἴσον)'. Euclid uses the same words in his definition VII.18: 'τετράγωνος ἀριθμός ἐστιν ὁ ἰσάκις ἴσος ἢ ὁ ὑπὸ δύο ἴσων ἀριθμῶν περιεχόμενος.' However, the usual modern translations do not take into account the ancient point of view: when it is translated as the 'product of two integers', for instance $m$ and $n$, the reference is that $p = m \times n$, and its meaning is to be a two variable ($m$ and $n$) operation. But the Greek wording used by Theaetetus, in agreement with Euclid's definition of 'to multiply' ('πολλαπλασιάζειν') in definition VII.15, [4] means it is an operation depending of $m$ and performed on $n$: to add $m$-time $n$ to itself. [5] (cf. also note 7), *infra*). A second point is that Theaetetus does not speak of the integers that **are** 'an equal times an equal', but more carefully that have 'the **power** to be' so. It is an important precaution, because any perfect square integer may also be written as the product of the unit by itself, [6] and if the 'power' was forgotten, any such integer will belong to the second class they define just afterwards, the 'rectangular integer' (cf. *infra*, points 7) and 9)). [7]

4) 'square' ('**τετράγωνόν**'). Dans le cadre des mathématiques grecques anciennes, 'τετράγωνός' signifie très généralement 'square'. Mais, *a priori*, le terme peut désigner plus généralement n'importe quel quadrilatère. On pourrait penser que 'equilateral' ('ἰσόπλευρον') sert à lever l'ambiguïté de τετράγωνόν, en partie au moins, car un losange, comme un carré, peut avoir quatre côtés ('πλευρά') égaux. [8] Mais ce n'est pas le cas, puisque le verbe utilisé un peu plus loin ('τετραγωνίζουσι') dont dérive 'τετράγωνός', renvoie bien la construction d'un carré associé à un rectangle. Le terme 'ἰσόπλευρον' est employé ici par symétrie avec 'προμήκη' qui signifie 'rectangle' pris au sens **strict** d'une figure rectangulaire dont les côtés sont **inégaux**. Dans le premier paragraphe du passage, les figures géométriques sont des rectangles (au sens large) dont les côtés sont égaux, donc des carrés. Dans le second, des rectangles (au sens large) de côtés inégaux. Dans les deux cas, il s'agit d'associer à un rectangle (au sens large) un carré. En fait, la double construction des jeunes garçons est d'une certaine façon redondante.

5) '**τὸν τοίνυν μεταξὺ τούτου**'. Le singulier pose un problème. Il faut considérer que 'τούτου' est un neutre indiquant l'intervalle laissé, dans le nombre, par la classe des carrés. C'est la seule occurrence de cette formule dans le corpus platonicien.

6) 'but is the product of a greater times a lesser, or of a lesser times a greater' ('**ἀλλ' ἢ πλείων ἐλαττονάκις ἢ ἐλάττων πλεονάκις γίγνεται**'): c'est une formule symétrique à celle définissant le carré au point 4), *supra*. Il s'agit des entiers qui ne sont pas des carrés, et donc ne **peuvent** être mis que sous la forme d'un produit d'un entier un nombre de fois différent de cet entier : soit un plus grand par un plus petit ou un plus petit par un plus grand. En écriture moderne, cela s'exprimerait par : les entiers de la forme $m \times n$, où $m$ et $n$ sont des entiers avec soit $m$ plus grand que $n$ i.e. $n$ plus petit que $m$, soit $m$ plus petit que $n$, i.e. $n$ plus grand que $m$. Pour un moderne, Théétète fait une longue périphrase pour dire que les deux nombres sont différents. Mais la définition de la multiplication des entiers comme additions successives nécessite ce double traitement (cf. *supra*, point 3)). Once again, the precaution to speak not about integers that are product of unequal integers, but that **cannot** be otherwise is noteworthy. Otherwise, as we remarked in point 3) above, any square integer can be a product of two unequal integers: the unit and itself.

7) 'side' ('**πλευρά**'): this term comes from geometry, meaning the side of a figure, for instance of a triangle or of a square. It is used by extension in arithmetic. According to Euclid's *Elements* book VII, whose language is very close to Theaetetus', an integer product of two integers is called a 'plane' ('ἐπίπεδος') integer, while each of these two integers in the product are its sides ('πλευρά') (definition 16). Using modern symbolism, let $p = m \times n$ be the product of the integers $m$ and $n$, then the 'sides' of $p$ are the two integers $m$ and $n$, and called 'factors' of $p$. [9]

However, the definition itself may entail a dangerous ambiguity. When $p$ is considered as the result of the product ($m$-time $n$), the sides are well-defined as the two integers $m$ and $n$. But it is also possible to understand that the sides are the sides of the integer itself, independently of the product. Then, at least for non-prime integers, they are not well defined for there are several different possibilities to write such an integer as a product of two integers; for instance, 15 is equal to 3 times 5 but also to 15 times 1 (the other possibilities as 5 times 3 or 1 times 15 give the same geometrical figure, so that the ambiguity is less of a problem). An ambiguity which seems to be at the root of the faulty proofs in book VIII of Euclid's *Elements* (cf. *infra*, note 12).

8) '**μείζων δὲ καὶ ἐλάττων ἀεὶ πλευρὰ αὐτὸν περιλαμβάνει**': Notre traduction rend une phrase active au passif pour des raisons de clarté. La traduction littérale serait : 'a greater side and a smaller one always encompasse it'.

9) 'rectangular shape' ('**προμήκει αὖ σχήματι**'): this shape is opposed to the square figure, so that, 'rectangular' has to be understood according to a strict meaning, as coming from a 'non-square rectangle' i.e. whose sides are unequal. [10] The plane figure representing the product of two integers in Euclid's *Elements*, book VII, is a rectangle (cf. point 4), *supra*) but with a larger definition than the one given by the boys, for it includes the squares. Euclid's definition avoids, as noted at the end of point 4), the redundancy in the boys' construction. [11]

10) 'a rectangular integer' ('**προμήκη ἀριθμὸν**'): [12] same remark as the above. In modern language, the 'rectangular integer' can be written **only** as a product of unequal integers. As noted above, Euclid does not consider separately the 'rectangular integers' in the strict sense of Theaetetus and the perfect squares, for they are both included in his definition of a 'plane integer' (cf. point 7), *supra*), though he gives a separate definition (VII.18) for the latter, almost identical to Theaetetus' (cf. *supra*, point 3)).

11) 'lines' ('**γραμμαί**') in Greek mathematics mean generally 'finite segments of straight lines', not as in modern mathematics, 'infinite lines'.

12) 'the rectangular' ('**τὸν ἑτερομήκη**'). Theaetetus uses it for 'rectangular' integer called previously '**προμήκη ἀριθμὸν**'. [13]

13) 'all the lines squaring the equilateral and plane integer' ("**Ὅσαι μὲν γραμμαὶ τὸν ἰσόπλευρον καὶ ἐπίπεδον ἀριθμὸν τετραγωνίζουσι**"). This means the sides of the square associated to any perfect square integer. Similarly, 'the [lines] squaring the rectangular [integer]' are the sides of the square of the same areas as the (strict) rectangles associated to the 'rectangular integers'. For instance,
   - To the square integer $9 = 3^2$ is associated the geometrical square of side 3, thus of area 9.
   - to the 'rectangular integer' 3 is associated a rectangle of sides 1 and 3, and the lines 'squaring three' is the side of the square of area 3, or in modern language the square root of 3.

14) 'not commensurable as length with the formers, they are commensurable as areas' ('**οὐ συμμέτρους ἐκείναις, τοῖς δ'ἐπιπέδοις ἃ δύνανται**'). The meaning is that the sides of the squares squaring the 'rectangular integer' are not commensurable to the sides of the square associated to a 'square integer', but (the areas of) all the squares are commensurable with each other. In modern terms, for instance the square root of 3 is not commensurable to any integer, but its square i.e. 3 is an integer, thus commensurable to any integer. [14]

15) 'And it is the same in the case of solids.' ('**καὶ περὶ τὰ στερεὰ ἄλλο τοιοῦτον.**'). What is the 'same' thing that Theaetetus is speaking about? To get an answer, we have to consider how the boys proceed in the former case, replacing the plane squares by the spatial cubes, the plane rectangles by the parallelepipeds, and the sides of the plane figures by the sides of the spatial ones. So we get something like this:
   i) The beginning is the same: the division of 'the integer in its totality in two parts'.
   ii) Then on the one hand, we have to separate instead of the integers that are 'product of an equal by an equal', the integers that are 'product of an equal by an equal by an equal' i.e. that are perfect cubes (in modern notations of the form $m^3$, where $m$ is an integer) from the non-perfect cubes. As a matter of fact, it is the definition in Euclid's *Elements*. [15]

iii) Each of these integers will 'be likened to the cube as a figure, and named cube and equilateral'.

iv) On the other hand, we get all the other integers, i.e. the ones which have not the power to be of the former kind, in other words, all the integers which are not perfect cubes.

v) To any integer of the first class i.e. the perfect cubes, we associate the geometrical figure of the cube of volume this integer. Thus to $n = m^3$, we associate the cube of side the integer $m$.

vi) The other integers, the ones 'in between as three and five' i.e. the non perfect-cubes, can only be the product of at least two different integers, thus 'encompassed by at least two different sides', or in modern language 'two different factors'. Then they are likened to a parallelepiped which has at least two different sides. And they will be called 'parallelepipedal integers'.

vii) To any such 'parallelepidal' integer will be associated a geometrical spatial figure, a parallelepiped of sides one, one, and this integer. [16] Thus, the volume defined by this spatial figure is this integer.

viii) The lines 'cubing' the first class of integers (the perfect cubes) which are the sides of the figure associated to them should be called 'length', as in the case of the squares.

ix) The lines 'cubing' any integer in the second class are the sides of the cube equal to the parallelepiped likened to this integer. This entails that, at the beginning of the 4$^{th}$ century, young boys as Theaetetus and Socrates had no problem to construct a cube of any volume. Moreover, while these sides are not commensurable 'as length' to the first ones (defined in v), *supra*), they are commensurable 'as the cubes they have the power to produce'. [17]

x) Thus, once again, it would make sense to name 'length' the first class, and 'powers' the second class. As a matter of fact, the 'length' is the same as in the case of the squares, since it gathers together all the integers. [18] And while these new spatial 'powers' are different from the plane ones (as for instance $\sqrt[3]{3}$ and $\sqrt{3}$), the reasons given by the boys to call them such are the same.

xi) Hence, for Theaetetus, 'powers' have a larger meaning than in the usual interpretations, and that it refers not only to the sides of the squares corresponding to non-square integers, but also to the sides of the cubes corresponding to non-cube integers. [19]

xii) From point ix), we get another consequence: in 399, the generalized solution of doubling the cube was enough well-established so to be known by young boys learning some elementary mathematics. This is an important testimony in any attempt to date the solution of the problem of doubling the cube, assigned usually to Archytas of Tarentum.

16) 'Excellent my boys' ("Ἄριστά γ' ἀνθρώπων"). En grec ancien, c'est un éloge appuyé ; chez les hommes on ne peut faire mieux. Mais chez Platon l'éloge est ambigu. On ne peut faire mieux chez les hommes, mais chez les dieux oui. Or le philosophe doit tenter de s'assimiler à dieu. En fait, Socrate indique ici que Théétète est loin du compte. L'erreur mathématique montre qu'il est nécessaire pour que Théétète engendre une théorie cohérente de l'irrationalité quadratique, il soit délivré des mathématiques à la Théodore. Pour accomplir ce qu'il a décrit dans sa pseudo-définition des 'puissances', il doit abandonner une mathématique fondée sur les dessins, pour une mathématique fondée sur les raisonnements. C'est ce type de mathématique qui permettra d'obtenir une théorie globale des racines quadratiques d'entiers, et plus généralement d'une large quantité de grandeurs irrationnelles, fondée non sur des algorithmes et des graphismes, mais sur des propositions analogues à celles obtenues au livre VII des *Éléments* d'Euclide (cf. *infra*, **Annexe**). La maïeutique socratique ouvre la voie à cette innovation, due à Théétète, en le délivrant d'une conception de la science fondée sur les choses sensibles, implantée par Théodore qui l'a lui-même reprise de son ami Protagoras (cf. *infra*, 0.1).

### 3. The meaning of '*arithmos*' ('integer') and '*dynameis*' ('powers') revisited

For the modern reader, the passage contains several traps, not least the one concerning the term '*arithmos*'.

**1) *Arithmos* (integer)**

Its almost universal translation by 'number' in English, '*nombre*' in French, '*Zahl*' in German, '*numeri*' in Italian, etc. makes sense, because 'arithmetic' means actually 'theory of numbers' (in French '*théorie des nombres*', in German 'Zalhentheorie', '*Teoria dei numeri*' in Italian, etc.). On the other hand, the term '*arithmos*' is given by Euclid as the second definition of book VII of the *Elements*: a 'collection of units'. It was already defined previously in almost the same terms by Aristotle. [20]

Although we usually find in modern translations a warning that, a 'number' was necessarily a 'positive non-zero integer' in the Antiquity, this warning is inadequate, and anyway not taken into account seriously. It is all the more inadequate for the analysis of this part of the *Theaetetus*, for the different kinds of 'number' (in the modern sense) are at the center of the problem.

i) **The usual translation of '*arithmos*'**

In modern mathematics there are many different notions of numbers, so that it is nowadays an ambiguous term. Its meanings include 'rational number', 'irrational number', 'real number', 'complex or quaternion number', though its usual meaning is 'real number'. In other words, in modern mathematics, the usual 'numbers' used to measure any length, weight, volume… are 'real numbers' e.g. 1.3, √2 or the constant π. They are used everywhere, even if their rigorous definition is far from easy.

Instead of such 'numbers' that were not defined until the end of the 19th century CE, [21] ancient Greek mathematicians used straight lines to represent them. One example is evidently the aforementioned incommensurable ratios. While these ratios cannot be represented within the framework of the integers, they can be represented within the framework of straight lines. For instance, as shown in Plato's *Meno* (84a, 85a-b), whereas it is extremely difficult to define the length of the side of the square doubling a square of two foot side, it is easy to obtain a geometric representation of it, by drawing the diagonal.

ii) **Our translation**

These problems explain, *inter alia*, why this part of the *Theaetetus* is difficult to interpret and have entailed numerous mistakes in the reading of the entire dialogue. Hence, we **will never use the term 'number',** with or without qualification, to designate a notion of ancient Greek mathematics.

## 2) In search of a definition of '*dynameis*' ('powers')

To move from an unlimited number of powers to one and only one definition [22] is a typical Socratic practice outlined already in Plato's early works. [23] The boys' definition of 'powers' is done in three steps. Among the integers, they consider, the perfect squares: 1, 4, 9, 16, …, the rectangular integers: 2, 3, 5, 6, 7, …, and finally, the relation between these two classes.

i) **The two classes of '*arithmos*'**

To understand Theodorus' lesson and to get further, Theaetetus and his young friend *Socrates* decide to review it together. Since Theodorus did not go further than 17, they tried to get a generalization to all integers. In a first step, they divided 'the integer in its totality' (cf. point 2.2), *supra*) in two disjoint classes, as follows: an integer $N$ belongs to either of these classes, depending on whether or not it 'has the power to be a product of an equal times an equal' i.e. (cf. point 2.3), *supra*) if there is an integer $k$ such that $N = k^2$.

In the first case, the integer is associated to the following geometrical figure: a square whose area is this integer. For instance, $N = k^2$ is associated to the square of area $N$ i.e. of side $k$. The boys called such an integer 'square and equilateral' (cf. point 2.4), *supra*). In the second case, they claim that, since any such integer $n$ can be written as a product of two unequal integers $p$ and $q$, it can be associated to a rectangle in a strict sense i.e. with unequal sides $p$ and $q$. If, because of its definition, it is clear that $n$ cannot be the product of two equal integers, it is not evident it needs to be a product of two integers. [24] However, when 1 is regarded as an integer, any integer may indeed be written as a product of two integers (1 and this integer), so that any non-square integer is actually a 'rectangular integer' ('προμήκη ἀριθμὸν' cf. point 2.6) and 2.10), *supra*') with one side longer/shorter than the other ('μείζων δὲ καὶ ἐλάττων ἀεὶ πλευρὰ').

All these 'rectangular integers' are located between the 'square' (*sic*, cf. 2.5)), and we have the following sequence, where the 'square' is written in bold type:
**1** 2 3 **4** 5 6 7 8 **9** 10 11 12 13 14 15 **16** 17 18 19 20 21 22 23 24 **25** …

This definition may entail an ambiguity. Indeed, there are several different ways for an integer to be represented as a product of two integers (cf. point 2.7), *supra*), thus as a geometrical rectangle. However, the problem disappears when we consider Theodorus' lesson, for the construction is carried out using only rectangles of which one side is equal to 1 foot (cf. B-O, Figures 4, 5, 6).

Let us also remark that the constructions need to be different according to the two classes. Certainly the construction for 'rectangular integers' may be extended to the 'square' integers, however it would be different from the construction given above for the latter. For instance, to 9 would be associated instead of a square of side 3, a rectangle of sides 1 and 9. [25]

ii) **The meaning(s) of '*dynameis*'**

At the beginning of the mathematical part, Theaetetus speaks of 'powers'. The meaning of this term has triggered a long controversy among modern scholars: does it means 'squares' or 'sides of squares'? Its first use in the passage (cf. B-O, point III.2) is in a vague sentence, introduction to Theodorus' lesson ('Περὶ δυνάμεων', 147d1 and implicitly in 147d2 ff '…τρίποδος πέρι καὶ πεντέποδος…'). It is not really conclusive, but since it concerns the construction of a figure, it would better correspond to the knowledge of the sides. Its second use, when Theaetetus speaks of areas ('three', 'five', …, 'seventeen'), is inconclusive as well, since it may refer either to the surfaces themselves or their sides.

In the passage of *Theaetetus* (147d7-148b2) studied here, the word appears several times, beginning with a remark that 'the powers are unlimited in plurality'. It is true for squares as well as for their sides, thus the remark apparently provides no argument for either branch. Nevertheless, let us consider it more in detail.

- If the 'powers' were geometrical squares, the sentence would mean that all the multiples of the unit-square are quantitatively unlimited, in other words that the integers are unlimited. This is so trivial that any explanation would be superfluous and that the verb '*phanein*' would not be used here.

- Thus, the property of being 'quantitatively unlimited' must be some clear consequence of Theodorus' lesson, but not a pure triviality. [26] According to their definition of 'powers' given later, [27] they understood from the lesson that all the square roots of non perfect-squares integers are incommensurable to the unit. Thus, Theaetetus' claim is equivalent to contending that there is an unlimited quantity of such integers i.e. non perfect square integers. [28] Though it is clear, it is not a pure tautology and it is worth a remark in passing. The conclusion is that 'powers' mean here sides of squares, but only these sides incommensurable to the unit. [29]

a) Then the boys begin a classification of the sides of squares into two classes (cf. §i), *supra*). One is composed of the sides of squares whose areas are perfect square integers, corresponds to the integers; [30] the other is composed of the sides of squares whose areas are not perfect squares so that they are 'not commensurable as length with the formers, [but] they are commensurable as areas that they have the power to produce' i.e. as areas of squares (148b1-2). In both cases 'powers' necessarily mean the sides of squares, not the squares themselves. In modern mathematical language, they are the irrational square roots of integers, e.g. √2, √3, √5, and so on. However, the boys claim to be defining 'powers'. Thus, there must be some changes between the beginning of the text (147a1-5) and this second part (147a6 ff.). It is highly unlikely that the boys either completely changed the usual meaning, or replicate it. Their purpose is to specify a seemingly overly general and vague sense of 'powers', which included the construction of the side of squares whose areas are integers. [31]

According to the boys, the 'powers' are the 'incommensurable sides' of the squares that they associated to integers.

Hence, on the one hand, they defined the 'length' as the sides of squares (148a7) represented by perfect square integers, thus these sides are integers, all the integers (cf. *supra*, note 30). On the other hand, they defined the 'powers' [32] as sides of squares of areas represented by non-perfect square integers. They claim that: i) the sides of these squares (in modern symbolism e.g. √2, √3, √5) are not commensurable with the formers, i.e. the newly defined 'length' (e.g. 1, 2, 3); ii) the squares built on these sides (the 'powers') are whole integers, and obviously commensurable with the squares built on the 'length'. [33]

    b)    To conclude, we may say that for Theaetetus and his friend *Socrates*, the meaning of 'powers' is definitely closer to 'sides of the squares' than to 'squares'. However, the term 'powers' does not maintain the same meaning everywhere. As claimed by Aristotle (cf. *supra*, note 31), it is linked to graphical representation of geometrical means. Since it is done by drawing a square, an ambiguity between area or side does not really matter, for the drawing of a square gives them altogether. When a distinction is needed, it is possible to add 'as length' (147d4) or conversely 'as area' (148b2), to make clear it is about the sides or the areas. [34] And in the last part, when Theaetetus gives the definition of 'powers', the word refers to the sides of the squares whose areas are a non-square integers.

As a matter of fact, the boys never defined what a 'power' **is** (the singular opposed to the plural 'powers'), as should be expected according to Theaetetus' analogy between their work and the philosopher's inquiry, a definition of 'what **is** knowledge'. [35] In both cases, they do not give a **definition** either of 'power' or 'knowledge, they merely give **exemplifications**.

Plato suggests that Theodorus' courses do not lead his young pupils, including the best and the brightest, to understand really mathematics. Moreover, as we will see, in the last part of the passage, Plato shows us Theaetetus and his friend *Socrates* making some mathematical mistake stemming from the graphical kind of mathematics taught by Theodorus.

### 4. The circularity of *Socrates*-Theaetetus' project

The boys called the side of some squares 'powers' for two reasons.

    a)    First, they are **incommensurable** with the 'length' (148b1) previously defined as the sides equal to integers (cf. *supra*, 2).ii).a)).

    b)    Second, the squares built on their sides are **commensurable** with the 'length' (i.e. the integers) for the areas of these squares are integers.

As explained in paragraph 3.2), the two classes of sides of squares correspond to the square roots of the two classes of integers defined previously, respectively the perfect square integers, and the non-perfect square integers. [36] In other words, the division is carried out between

- on the one hand, the class of integers (the 'length')
- on the other hand, the class of non-integer square roots of integers (the 'powers').

i) Since all the areas of the squares considered by the boys are integers, the point b) is a triviality.

ii) Moreover, it is implicitly assumed by the boys that these two classes are disjoint, a property entailed by the following trivial equivalence (written in modern terms):

the square root of an integer is an integer is equivalent to this integer is a perfect square. [37] This disjunction is also the consequence of the above point a), since all integers are commensurable to one another. So far so good.

iii) Now let us consider the above reason a) given by the boys to call 'powers' the objects of the second class. Translated in in modern terms, we get the following equivalence:

The square root of an integer is rational is equivalent to this integer is a perfect square. [38]

We will consider it in detail in the **Appendix**, and show that it is a highly non-trivial result. It is impossible to consider such a result so evident to be used in the framework of a definition.

Nevertheless, let us remark that the equivalences in ii) and iii) seem pretty close to one another, the only difference is that in iii) the term 'rational' replaces the second 'integer' in ii). No wonder that both have been often confused. As a matter of fact, Theaetetus' account of Theodorus' lesson (cf. B-O) suggests strongly that the boys made this confusion, for the lesson presents only two kinds of magnitudes: on the one hand the integers, on the other hand the irrational, without the middle term of rational magnitudes. [39] In other words, Plato presents Theaetetus stating a 'definition' that needs the use of the highly non trivial - including for modern mathematicians - '**equivalence number 2**' while thinking to the trivial - including for the Greeks of the 5th century - '**equivalence number 1**'. [40]

Following the way Theodorus taught mathematics, it can be easily understood why the boys made the mistake in their definition. However, we are then confronted with a puzzle: why did Plato present the young Theaetetus and his comrade *Socrates* as providing a (true) statement without having a hint about the necessity this statement needs a (difficult) proof?

An unsatisfactory answer given by some historians is that it is a way to pay tribute to Thaetetus by making him just saying a statement that he will prove some years later, once a mathematician; in other words, Plato makes Theaetetus prophesy. This claim goes directly against what Socrates says a little later, that mathematics is a science where everything has to be proved (162e4-163a1). Actually, the dialogue is far from praising Theodorus and Theaetetus: both are supporters of the sophist Protagoras; Theodorus' mathematics is described using the very graphical tools blamed by Socrates in the *Republic* (VII, 527a); the rest of the dialogue is nothing else than a sequence of 'apories', Theaetetus (and Theodorus) being unable to solve any problem raised by Socrates.

To get an answer, we have to look in a different direction. Let us first remark that the boys' classification, in line with Theodorus' lesson, starts from the 'integer in its totality', and they get two classes of sides (of squares): the 'length' i.e. the integers, and the 'powers' i.e. in modern language, the irrational square roots of integers. Such a division does not exhaust all the cases of commensurable sides (to the unit), contrary to the definitions X.2-3 of Euclid's *Elements*. [41] This awkwardness was already noticed in the Antiquity. [42] However, from Plato's point of view on division, it is not only awkward, but a mistake. [43] It explains why the boys are blind to the principal difficulty in the problem. And nevertheless, everything they claim is true. As a matter of fact, Plato presents us here with a kind of error difficult to identify: to believe some point so trivial that it does not need a proof, i.e. the reason why the facts are such and not otherwise (*Meno* 98a), the worst kind of error for Plato (*Apology* 22d-e; *Theaetetus* 150c-e; 210c). In any case, our analysis makes it clear that the solution was not within the reach of the boys (see *infra*, **Annex**).

The solution of the previous puzzle entails the answer to a more general one: why does Plato put the whole 'mathematical part' inside a philosophical inquiry, or in other words what is the **philosophical** purpose of this passage?

## II. The philosophical significance of the mathematical part

We will consider briefly in this second part, first the common basis of the most usual interpretations (the so-called MSI) about the theory of knowledge ('*epistêmê*'), and show it is not consistent with Plato's text. Then, we will give the interpretation entailed by the above analysis of the 'mathematical passage'. We will conclude with what can be learnt about 'knowledge' ('*epistêmê*'), the subject of Socrates' inquiry, from the dialogue.

### 1. The 'Main standard interpretation' and Plato's Forms

The goal of this currently widespread interpretation is to claim a definition of knowledge that does not take account of the Forms, or of the soul, considered as elements which are obsolete. This position is rooted in an initial certainty that by the time he wrote the *Theaetetus*, Plato had abandoned the doctrine of Forms that can only be apprehended by a soul capable of separating itself from the body. [44] This leads them to translate '*eide*' by 'kinds' and to overlook three important passages: the 'mathematical part', the one on maieutics (*Theaetetus* 148e-151d), and another one which, following Plato himself, is described as a 'digression' (*Theaetetus* 172c-177d) [45], where Forms have to be taken into account.

#### i) Forms and maieutic

Since we studied the 'mathematical part', let us consider the text on maieutic. Dans la MSI, on estime que ce qui est dit de la maïeutique dans le *Théétète* est purement rhétorique. C'est une grave erreur. Socrate affirme d'une part ne rien avoir à transmettre à Théétète, et d'autre part, par ses questions, il propose de purifier l'âme de Théétète des erreurs qu'elle recèle. The necessity for Theaetetus to go through Socrates' maieutic has been already emphasized in the first part. [46] To obtain what he claimed in his pseudo-definition of powers and get a global theory of quadratic roots of integers, he will have to abandon Theodorus' mathematics, founded on algorithms and drawings, for one close from what is found in book VII of Euclid's *Elements* (cf. *infra*, **Appendix**). As a matter of fact, the text on Socrates' maieutic begins immediately after the mathematical passage.

As a midwife of the soul, Socrates explains that he can bring his interlocutor to give birth to a truth or promote a miscarriage if a mistake is uttered. In the *Theaetetus*, three definitions of knowledge are rejected, one after the other, because they prove to be incoherent. Socrates limits himself to asking questions and making remarks, while Theaetetus is leading the investigation. The reason is that Socrates has no knowledge to transmit: all he can do is to make his partner [47] give birth to a truth, as Socrates' mother helped other women give birth. Yet refutation is only one of two faces of the same coin, the other one being recollection. According to Socrates, his partner does not receive a truth transmitted by a teacher, but he has to discover within and by himself a truth (150d6-8). In another dialogue with which the *Theaetetus* is closely connected, the *Meno*, where one finds also a mathematical passage. [48] Its subject is about the question of how can virtue be found. Virtue cannot be the result of a natural gift, of practice or of teaching; it is obtained through 'recollection' as a divine gift. It should be compared with a passage from the *Phaedo* (72e-77a) in which the Forms are designated as the objects of recollection that the soul has contemplated in the past. [49]

### ii) Forms and the 'Digression'

Let us now consider the 'digression'. The forms are implicit everywhere in the passage. For instance, how can one fail to think of the Forms, when one reads such expressions as: 'as long as they reach that which is' (172d9 and Burnyeat (1990)); 'whatever a man is' (ib., 174b3); 'to come to the examination of justice itself and injustice: what each of them is, and in what do they differ from everything else or from one another' (ib., 175c2)? How could one not admit that the soul exists, and subsists after death, when one reads 'But is not possible, Theodorus, that evil should be destroyed – for there must always be something opposed to the good; nor is it possible that it should have its seat in heaven. But it must inevitably haunt human life, and prowl about this earth. This is why a man should make all haste to escape from earth to heaven; and escape means becoming as like god as possible; and a man becomes like God when he becomes just and pious, with understanding'. (ib., 176a-b); or 'And if we tell him that, unless he is delivered from this 'ability' of his, when he dies the place that is pure of all evil will not receive him; that he will forever go on living in this world a life after his own likeness – a bad man tied to bad company' (ib., 176e-177a).

Everything indicates that these two last passages should be interpreted as Cornford did, in a sense that represents a tradition going back to Platonism in Antiquity. [50] It entails that the MSI's claim that the Forms contemplated by the soul during a previous existence do not appear in the background of the *Theaetetus* is nothing more than a *petitio principii* based on a supposed evidence: the abandonment by Plato of the hypothesis of the Forms as the base of his theory of knowledge, something found nowhere in his texts. It is all the more lacking any credibility when we consider that it is followed by the *Sophist*, in which it is impossible not to admit the existence of the Forms and of the soul. [51]

### iii) Geometry and knowledge of the forms in the Republic

Our interpretation finds a confirmation in the famous discussion between Socrates and Glaucon on geometry in the *Republic* (VII 527a1-b12):

S. – They give ridiculous accounts of it, though they can't help it, for they speak like practical men, and all their accounts refer to doing things. They talk of 'squaring' 'applying', 'adding' and the like (λέγουσιν τετραγωνίζειν τε καὶ παρατείνειν καὶ προστιθέναι καὶ πάντα οὕτω), whereas the entire subject is pursued for the sake of knowledge (πᾶν τὸ μάθημα γνώσεως ἕνεκα ἐπιτηδευόμενον).

G. – Absolutely.

S. – And mustn't we also agree on a further point ?

G. – What is that ?

S. – That their accounts are for the sake of knowing what always is, not what comes into being and passes away (Ὡς τοῦ ἀεὶ ὄντος γνώσεως, ἀλλ' οὐ τοῦ ποτέ τι γιγνομένου καὶ ἀπολλυμένου).

G. – That's easy to agree to, for geometry is knowledge of what always is (τοῦ γὰρ ἀεὶ ὄντος ἡ γεωμετρικὴ γνῶσίς ἐστιν).

S. – Then it draws the soul towards truth and produces philosophic thought by directing upwards what we now wrongly direct downwards (Ὁλκὸν ἄρα, ψυχῆς πρὸς ἀλήθειαν εἴη ἂν καὶ ἀπεργαστικὸν φιλοσόφου διανοίας πρὸς τὸ ἄνω σχεῖν ἃ νῦν κάτω οὐ δέον ἔχομεν).

G. – As far as anything possibly can. [52]

Already in the *Republic*, geometry is linked to knowledge dealing not with becoming, but with being, that is with the Forms. A new argument that the 'Main standard interpretation' is out of the pale.

## 2. Our interpretation

Let us now return to the mathematical part and the consequences of its analysis for the understanding of the whole dialogue. The defect that led *Theaetetus* to the three successive failed attempts of definition of knowledge is already found in the mathematical passage. Knowledge can neither be sense-perception of geometric figures, any more than the supposed definition of 'powers' can be a mathematical definition, though it respects all the logical criteria of truth. *A fortiori* it is not a scientific one, for it is not a definition at all.

No definition of knowledge will be given by Socrates, who affirms that he knows nothing. [53] As he explains in the famous passage on maieutics (148e-151d), Socrates, who cannot conceive, has received the gift of helping others to give birth. He is therefore attempting to deliver Theaetetus of the pains of childbirth concerning the definition of knowledge. This is not a slightly ridiculous *mise en scène*, intended to surprise or even amuse the reader. [54] Maieutics refers to recollection, which in fact implies an original conception of education. Plato illustrates this in the *Meno*, by questioning an illiterate servant regarding the duplication of a square: to learn is not to receive an item of information transmitted by a teacher, [55] but to recollect an item of knowledge previously acquired by a soul that existed prior to its entry into such-and-such a body. In the *Theaetetus* Socrates tries to help a young boy give birth to a definition of knowledge, and the dead end of the attempted definitions was anticipated in the definition of powers in the mathematical part.

Theaetetus makes three successive attempts of very unequal lengths.

i) Knowledge is perception (151d-186e) = 35 Stephanus pages

ii) Knowledge is true opinion/judgment (187a-201c) = 14 Stephanus pages

iii) Knowledge is true opinion/judgment accompanied by a *logos* (201c-210b) = 9 Stephanus pages

a) The *logos* is the vocal expression of true opinion/judgment (206d-e)

b) The *logos* is the complete enumeration of the elements (206e-208c)

c) The *logos* is the designation of the characteristic difference (208c-210b)

All these attempts, however, collapse under Socrates' examination.

## 3. The meaning of 'knowledge' ('*epistêmê*')

### i) Knowledge is sense-perception ('*aisthêsis*') [56] (151d-186e)

Theaetetus' first definition of knowledge is 'perception/sensation' ('*aisthêsis*') (151e3). To a modern, this seems, *prima facie*, to be a strange answer coming from a mathematician. However, it is consistent with Theodorus' practice, working and providing proofs through drawings. Indeed, as previously seen, he congratulated himself explicitly for having, unlike Socrates, left 'bare discourses' for geometry as soon as possible (165a). His method consists in making the results appear to his students through figures. Its advantage is the ease with which it enables that solutions of problems and proofs appear almost as evident. [57] There is however a drawback: its primitive aspect as compared to the mathematics carried out just a decade later. The latter includes a theory of irrationals founded on arithmetical reasoning, ascribed to none other than Theaetetus, impossible to obtain by graphic means. Another pitfall of Theodorus' graphical geometry is that it entails that not-being is not an object of mathematics, for only what exists can be represented graphically. This very conception led *Socrates*-Theaetetus to consider as a definition of 'powers' what is actually an extremely difficult proposition, requiring a long demonstration not available at the time. [58]

The same error that 'non-being' cannot be the object of the *logos* (for it is '*alogon*'), or that it is of no interest (for it is impossible to express, or '*arrêton*') is condemned, along a more general framework, in the *Sophist* (238e), a dialogue that is staged on the next day.

This definition of knowledge through sense-perception is placed under the patronage of the sophist Protagoras (151e-153c). In fact, Protagoras was merely formulating the position of other so-called 'thinkers of nature' (except Parmenides), particularly Heraclitus and Empedocles, and of such poets as Epicharmus, Homer, and Hesiod: nothing is, but everything is becoming (152a-155c). More subtly, it can be pointed out there is nothing stable about the action-passion relation, implied by the perception/sensation pair (155c-157c). The exposition ends by recalling three formulas (157c-160e): Heraclitus' universal flux, Protagoras' man as measure of all things, and Theaetetus' perception as knowledge.

Socrates sets forth several arguments against this. It should be noted that the first argument is addressed directly to Theodorus. If knowledge is sense/perception, all human beings will be equal (160e-163a), and knowledge will last no longer than perception (163b-165e). In fact – this is the second objection – Protagoras, whom Socrates defends just long enough to expose Protagoras' doctrine (166a-169d), contradicts himself when he admits that all human beings are not equal (169d-172c); in addition, this doctrine makes it impossible to speak of the future (177c-179c). The third objection (179c-184b) is radical: if everything is always changing, then there will be no longer any stable objects, and consequently no knowledge or discourse. The final attempt (184b-186e), which brings up maieutics, involves the soul. Since each sense provides only specialized and fragmentary information, an agent is required to unite them, give them meaning, and allow them to express themselves. This agent is the soul, which works directly on realities (187a), a paradoxical consequence of Protagoras' theory.

ii) **Knowledge is true opinion/ judgement ('*alêthês doxa*',187a-201c)**

Theaetetus' second attempt is: 'true opinion/judgment' (187b5) which both the finality and the end of the process of perception/sensation. [59] Unlike his first answer, this one is not spontaneous, but the result of Socrates' successive refutations. More precisely, once the identity between 'perception/sensation' and 'knowledge' has been destroyed, the boy is led to consider the soul through its 'powers', understood as faculties [60] (184e). Socrates advises him to make a new start ('obliterate everything said previously' [61]). Theaetetus then gives a new answer: knowledge is 'opinion/judgment ('*doxa*'), [62] going on to add the qualification of 'true' for there are also 'false opinions/judgments'. Socrates discusses Theaetetus' new definition at length, namely its addition of truth, wondering whether 'false opinion/judgment' would indeed be possible under this definition of knowledge. Then he shows briefly the complete failure of such a definition, through the whole part of an art despised by Plato, the art of the courts (for instance, *Gorgias* 464b-466b; *Phaedrus* 260a-262a).

In order to understand the discussion, we must specify the meaning of opinion/judgment ('*doxa*'). When the physical affections transmitted perceived as an intuition through sensation have reached the rational part of the soul (*Timaeus* 64b3-6), they triggered judgment, for instance: 'This apple is red' or in the *Sophist* (263a): 'Theaetetus sits' or 'Theaetetus flies', which is true or false. But because the objects of all knowledge of what is sensible are becoming, and hence always changing, 'opinion/judgment' pertains to conjecture. Opinion/judgment is thus an approximate thinking which passes itself off as a judgment (this is why '*doxa*' can also be translated by 'judgment') on what seems to be the case (the noun '*doxa*' comes from the verb '*dokein*', 'to seem, to appear'), or sensible things, which are mere images of the intelligible. It is liable to be true or false without ever being able to account for its truth or falsity and is therefore a mode of thinking intermediary between knowledge and ignorance. According to Plato, the only genuine knowledge will be that which deals with intelligible forms whose sensible particulars are but images. Therefore, it is only through knowledge of the intelligible that knowledge of sensible particulars will be possible. Hence, the object known by the intellect (the Form of horse) can bear the same name as the sensible thing (this horse), although the former and the latter are two distinct realities (*Phédon* 95e-102a ; *Republic*, VI, 509d-511a ; *Timaeus* 27d-28a among others).

In fact, this section raises the problem of error: how can a false opinion/judgment arise (187d-189b)? Error cannot be defined by substituting one object for another, for supposing that there is debate on the subject, opinion/judgment will not be in a position to decide which is the right one (189b-191c). For the same reason, error cannot be defined by substituting one memory for another, according to the image of a block of wax containing imprints, for how could one choose? (191c-197a). Finally, one must distinguish between having knowledge and possessing it. One may possess a large number of opinions/judgments in one's soul, like birds in an aviary, but when one wishes to have knowledge, it is impossible to know, by catching one in mid-flight, that one has caught the right one (197a-200d). A criterion is needed to decide which is the right opinion/judgment. However, even if such a criterion existed, it would require another one to decide whether an opinion/judgment verifies it; and then, again, another one to decide if the second criterion is verified by this opinion/judgment, and so on to infinity (200a-c). In short, opinion/judgment can never reach rigorous knowledge. Hence, it is simply impossible for even a true opinion/judgment to become knowledge.

The practice of the courts (199d-201c) shows decisively that it is impossible to differentiate between true and false through opinion/judgment. [63] A judge can be persuaded by a likely speech without knowing if it is true or false, for what result from the claims of two opposite parties may be a 'true opinion/judgment', but certainly not 'knowledge'. In the *Phaedrus*, one finds this practical example given by Tisias (5[th] century BCE) in his treatise on rhetoric: '…if a weak and punky man is taken to court because he beat up a strong but cowardly one and stole his cloak or something else, neither one should tell the truth. 'The coward must say that the spunky man didn't beat him up all by himself, while the latter must rebut this by saying that only the two of them were there, and fall back on that well-worn plea, 'How could a man like me attack a man like him ? The strong man, naturally, will not admit his cowardice, but will try to invent some other lie, and may thus give his opponent the chance to refute him. And in other cases, speaking as the art dictates will take similar forms.' (273b-c). [64] In a law court nobody cares about the truth ; one must say not what has actually happened, but what is likely in order to persuade the judges.

iii) **Knowledge is true opinion/judgment accompanied by logos (201c-210b)**

Let us now consider the third and last definition, for it differs from the other two by several features. It is not Theaetetus' own definition, but something he heard somewhere, told by someone, so that nothing is known about its origins. This aspect is emphasized further, for the boy, who is praised for his sharp memory ('*mnêmones*', 144a) claims he has nothing other than a dim recollection of it (201d7-8). For the first time the term '*logos*' appears out of nowhere (201c9-201d1). While Socrates' examinations of the first two definitions were long and painstaking, this one is very brief.

This last definition is introduced within the framework of two dreams. Socrates considers the new definition reported by Theaetetus to be a dream. Then, he proposes an exchange of dreams to the boy, Socrates' dream for the boy's dream. Both are recognized as 'dreams' because their authors are neither Theaetetus himself nor Socrates, since they became known to them by hearsay. The former heard it told by some unidentified speaker, the latter claims he believes he heard it from some unspecified people. The term used by Socrates ('ὄναρ ἀντὶ ὀνείρατος') is the same used for the state of the servant-boy at the end of the search for the line doubling a 2-foot side square (85c9). At this point, as Socrates states clearly to Meno, the servant-boy is ready for a full understanding of the solution, though he has not yet reached it and it still requires much work.

An obvious difficulty with the third definition is the meaning of '*logos*', a polysemous term that may signify explanation, discourse, reason (as a faculty), and causal reason, as well as mathematical proportion.

To Theaetetus' dream, Socrates adds an important point: the consideration of the logos as a whole (201c-203a), and as such, it must be broken down into its components. However, in order to consider these components, they have to be considered in turn as wholes, whose components must once again be determined, and so on to infinity. To avoid such an infinite regress, some components need to be primary elements given only by perception, and only having a name. Thus, though such components may be recognizable, they are not susceptible of reason, for reason, like definitions, implies a whole made up of knowable elements, as shown in grammar, in music and in the other sciences (203a-206c). Thus, this approach to the *logos* must be abandoned. And once again, the root of the question is found in the 'mathematical part'. As seen previously, [65] *Socrates*-Theaetetus' attempt of definition begins by the decomposition of any integer into a product of two factors. [66] Hence, if primary elements were unknowable, we would get the absurdity that it would be possible to know the product without knowing anything on its factors. [67] However, if such conception of the *logos* is wrong, how else can it be understood? Here, three answers are proposed successively, then rejected, by Socrates.

First, the *logos* cannot merely be the vocal expression of true opinion/judgment, since it would suffice to know to speak in order to get knowledge (206d-e). Nor can it be the exact and complete enumeration of its elements, for how could one know that such and such an object had them all, and moreover how could one be sure that one could define all these elements? The conclusion is that it is impossible to think of a compound either as only the sum of its elements, or as totally separated from them (204d-e); moreover, the same reasoning shows that a compound cannot be unknowable (205d-e).

Finally, defining the *logos* by designating the specific difference merely duplicates true opinion/judgment, which can only be true if it has designated the specific difference.

Consequently, none of the definitions of which Theaetetus will be delivered is viable. [68] Going back to the end of the 'mathematical part', we find a radical objection to all the last three definitions of 'knowledge': the downgrading of the 'truth' in favor of the '*logos*' as far as knowledge is concerned.

The boys' reasoning seems irreproachable from the logical point of view, with a true conclusion, thus, *a fortiori* true from the point of view of opinion/judgment. Nevertheless, the reasoning and the conclusion as a whole are incorrect because of the lack of a proof. Without a (correct) proof, *Socrates*-Theaetetus' supposed definition is nothing more than a ghost. [69] It has no more reality than the definitions of knowledge of which Socrates delivers Theaetetus through his art of 'maieutic', all of which appear to be mere smokes and mirrors (210a6-9). The boys' pseudo-definition of powers showed they were just on the border separating Theodorus' old graphical geometry from a new, more abstract geometry leading to the general theory of incommensurability of which only a small part will be concerned by Theaetetus' 'powers'. Similarly, at the end of the dialogue Socrates claims that Theaetetus is ready to inquire into knowledge. However, the question now seems to be not so much the search for its definition as for the definition of '*logos*' and its connection to 'true' and 'opinion/judgment'.

According to this conclusion, the sequel of the *Sophist* and the demise of 'absolute non-being' will be one step toward defining the '*logos*' in connection to 'knowledge'. [70] This is why Socrates, who speaks again of maieutics at the end of the dialogue (210b-d), invites Theaetetus to a discussion the next day. This discussion will be depicted in the *Sophist*, in which *Theaetetus* will play an important role, where the essential topic of discussion will be the *logos*. Between the two dialogues, a fundamental step has been taken: whereas in the *Theaetetus*, the problem studied was 'how can opinion/judgment be false?' (187d-200d), in the *Sophist* the goal is to answer the question 'how can a *logos* be false?' Despite the long and incredibly meticulous analysis by Socrates of the former question in the *Theaetetus* (192a-194a) no answer is found to the first question. Conversely, the *Sophist* gives one to the latter, but it presupposes the crisis of the 'being', and paradoxically the existence of the 'non-being'.

    For Plato, a definition of knowledge cannot be obtained by an investigation at the level of the opinion/judgment. Actually, the analysis of the mathematical part of the *Theaetetus* leads to a deep philosophical conclusion. Like the difference between rational and irrational magnitudes, the difference between opinion and knowledge are insuperable. The objects of opinion/judgment on the one hand and of knowledge on the other are not at the same level, so that nothing added to the former would transform it into the latter. And the definition 'true opinion and logos' is nothing else than an oxymoron. This conclusion clashes definitively with the prejudices of logical empiricism, but this is what Plato shows, with meticulous care, throughout the whole dialogue, and in a particularly evident form in its 'mathematical part'.

### III. Appendix. The mistakes in *Socrates*-Theaetetus' 'definition'

#### 1. The two problems

The text analyzed here has been generally interpreted, in the self-called 'Main Standard Interpretation, on the basis of two hypotheses. First, Plato is giving a tribute to both Theodorus and Theaetetus [71]; we have tried to show that it is not so obvious. [72] The second presupposition is about the boys' 'definition' of 'powers' (especially lines 148b1-2); there are several mathematical difficulties in their method, the most important being that they call 'definition' a (difficult) proposition. [73] It is all the more important to get a clear view of their 'definition' that it is essential for the interpretation of the puzzling end of the dialogue, the last definition of knowledge (201c-210b). [74] Among historians of mathematics, the standard interpretation of this pseudo-definition, is that it is a particular case of the generalized proposition 9 in book X of Euclid's *Elements*.

Let us recall the mathematical problem behind the boys' definition, before considering proposition X.9 of Euclid's *Elements* and comparing it to the proposition behind the boys' pseudo-definition.

i) As we tried to show previously, the problem is connected to the boys' (as well as some commentators') mistaken view that it needs only the following equivalence:

for some integer $k$, we have: $\sqrt{n} = k \Leftrightarrow n$ is a perfect square.

Since the second part of the equivalence is obtained simply by squaring the first (cf. *supra*, §I.4), it is a trivial equivalence ('**equivalence number 1**').

ii) The boys' supposed 'definition' is actually a proposition that requires the modern so-called 'Gauss' lemma'. [75] According to a long tradition, this 'definition' was indeed proved in Antiquity, but several years later, by Theaetetus alone, once a fully grown mathematician, without the help of his comrade *Socrates*. [76] As noted in paragraph I.4, *supra*, the problem lies in the claim that the 'powers are non-commensurable with [the 'length']'. [77] In modern mathematics, it means that the square root of an integer $n$ is rational if and only if the integer is a perfect square, [78] i.e.:

for some integers $m, p, q$, we have: $\sqrt{n}/m = p/q \Leftrightarrow$ for some integer $k$, we have: $n = k^2$.

('**equivalence number 2**').

Since no extant texts exist about the way Theaetetus is said to have obtained it, the only ones we can use are from Euclid's arithmetical books (VII to IX) and book X about irrational magnitudes.

#### 2. Proposition X.9 in Euclid's *Elements*

Let us consider the 'shocking unanimity' of some general 'strange delusion' among commentators and especially modern historians of sciences, [79] that Theaetetus' pseudo-definition is generalized in proposition X.9. [80]

The so called evidence is a text in a commentary attributed to Pappus of book X of Euclid's *Elements*. [81] The author compares Euclid's proposition X.9 to Theaetetus' work as accounted in the eponymous dialogue. He expresses his admiration for Euclid's proposition presented as much more general than Theaetetus'. Moreover, an anonymous scholium to proposition X.9 claims that this very proposition was discovered by Theaetetus, though in a less general form than Euclid's proposition. [82] For brevity and clarity, we will use the modern language, though it would be easy to translate it in a form similar of the one in the *Elements*. Proposition X.9 states:

'The squares on straight lines commensurable in length have to one another the ratio which a square number has to a square number; and squares which have to one another the ratio which a square number has to a square number will also have their sides commensurable in length. But the squares on straight lines incommensurable in length have not to one another the ratio which a square number has to a square number; and squares which have not to one another the ratio which a square number has to a square number will not have their sides commensurable in length either.' [83]

In modern symbolic language, the proposition may be translated as:

Let $c$ and $d$ be two lengths, we have the following equivalences:

$c/d = m/n$ for some integers $m$ and $n$ is equivalent to $c^2/d^2 = p^2/q^2$ for some integers $p$ and $q$.

Taking respectively for $c$ and $d$ the side of the squares of area $a$ and $b$, this can be written as:

$\sqrt{a}/\sqrt{b} = m/n$ is equivalent to $a/b = p^2/q^2$.

Just as the '**equivalence number 1**', this equivalence ('**equivalence number 3**') [84] is trivial, since the second part is simply the squaring of the first i.e. we can take $m = p$ and $n = q$. [85] Even if. in the context of Greek mathematics, this proposition is far from being self evident (cf. note 85 and 84), from the point of view of modern mathematics, proposition X.9 is a trivial one.

### 3. Comparison of the three equivalences
   a) According to '**equivalence number 2**':

for some integers $m, p, q$, we have: $\sqrt{n}/m = p/q \Leftrightarrow$ for some integer $k$, we have: $n = k^2$.

Then, let us remark the obvious analogy between this equivalence and 'equivalence number 3', since written in modern terms,

'equivalence number 2' claims:

- The square root of an integer is rational if and only if this integer is a square (of an integer), [86]

while 'equivalence number 3' claims:

- The square root of a rational is rational if and only if this rational is a square (of a rational). [87]

Since an integer is a particular case of (in the modern sense) a rational 'number', it seems *a priori* natural to think that the second statement is a generalization (to rational 'numbers') of the first one.

However, it appears impossible to prove '**equivalence number 2**' without using, in one way or another, the results of the theory of 'relatively prime integers' of book VII. [88] But this theory was not available at the time of Theodorus, since otherwise he would have been able to prove at once all the cases (as well as many others) he considered 'each case in turn'. [89]

Now, let us try to get the simplest proof for '**equivalence number 2**'. First, to summarize propositions VII.20-22 in one statement, [90] we get this **fundamental result**: For any ratio of integers $a/b$

- there exists two integers $c$, $d$ prime to one another such that $a/b = c/d$.
- Moreover, $c$ (respectively $d$) divides $a$ (respectively $b$).

Let us now return to '**equivalence number 2**'.

The equality $\sqrt{n}/m = p/q$ entails $\sqrt{n}/1 = t/q$ with $t = mp$. From the first part of the above fundamental result, there exist two integers $r$ and $s$ prime to one another such that $t/q = r/s$. By squaring both sides of the equality, we get:

$(\sqrt{n}/1)^2 =$ (by definition) $(r/s)^2 = n/1 = r^2/s^2$; from proposition VII.27, $r^2$ and $s^2$ are still prime to one another; thus from the second part of the fundamental result, $s^2$ divides $1$, thus $s^2 = 1$. Hence $n/1 = r^2/1$ i.e. $n = r^2$, and $n$ is a perfect square. Since the converse is evident, '**equivalence number 2**' is proved. [91]

    b)    Let us now compare '**equivalence number 2**' to the two other equivalences. While '**equivalences number 1** and **3**' are mathematical trivialities, '**equivalence number 2**' requires the theory of relatively prime integers i.e. a large part of book VII of the *Elements*. Hence, the former ('**equivalences number 1** and **3**') are certainly not identical with, and are even less a generalization of, the latter ('**equivalence number 2**').

    c)    All what had been said concerning the problems of the boys' pseudo-definition of the 'powers' as sides of squares can be extended, almost word for word, to a similar definition for the sides of cubes, since proposition VII.27 gives also the corresponding result for cubes. In modern language, exactly as in the proof of **equivalence number 2**, we get the following:

for some integers $m$, $p$, $q$, we have: $\sqrt[3]{n}/m = p/q \Leftrightarrow$ for some integer $k$, we have: $n = k^3$,

using the second part of proposition VII.27 stating that 'the cube of two integers are prime to one another if and only if these integers are prime to one another'.

Now, as we tried to show, the boys were thinking to 'equivalence number 1' whose extension to the case of cubes would be written as:

$$\sqrt[3]{n}/m = k \Leftrightarrow n = k^3,$$

Thus, once again, the parallelism is so evident and trivial, that the boys could consider its extension does not need any explanation (cf. our reconstruction of the cubic case, *supra*, §I.2.15), in particular points x) et xi)).

### 4. Conclusive remark

The conclusion is stunning: Pappus, a mathematician of the 4[th] century CE did no longer understand the mathematics of the 5[th] or the 4[th] century BCE, or at least the one described in Plato's texts. Since lines corresponding to an integer are a particular case of rational lines, proposition X.9 could be reasonably understood as a generalization of Theaetetus's claim about the formers. Someone like Pappus, always ready to praise Euclid, had no reason to look deeper, keen to show that Euclid obtained much better and general results than earlier mathematicians probably known through compilations or commentaries on Plato. This is supported by Pappus' testimony itself. In the introduction to his comparison, he speaks about 'those who have been influenced by speculation concerning knowledge of Plato', we can safely infer, Platonists or commentators of his dialogue. They, Pappus adds, 'suppose that the definition of straight lines commensurable in length and square and commensurable in square only which he gives in his book entitled, *Theaetetus*, **does not at all correspond** with what Euclid proves concerning these lines' (§10, Thompson (1930), p. 72). Thus, at the time of Pappus, there was still a strong tradition, at least among Plato's followers, that maintained the right consideration pointing to the difference between what is said in the dialogue and proposition X.9, but the mathematics reasons were probably lost.

More surprising is that modern historians followed Pappus without checking his claims beforehand.

l'irrationalité', *Bulletin de l'Académie Royale des sciences et des Lettres du Danemark*, année 1910, 5, p. 395-435 = Zeuthen (1910).

Notes

---

[1] The boy, friend of Theaetetus, namesake of the philosopher. To differentiate him from the latter, we will write his name in italics.

[2] The lesson was about some incommensurable magnitudes i.e. irrational quantities. For all the questions about commensurability/incommensurability and rationality/irrationality we refer to B-O, **Appendix III.**

[3] Though for the clarity, we will often use modern mathematical language to give some explanations about the mathematics involved by the text.

[4] 'An integer is said to multiply an integer when that which is multiplied is added to itself as many times as there are units in the other, and thus some integer is produced.' (All our translations from Euclid's *Elements* are from Heath (1908), sometimes slightly modified).

[5] Though the result is evidently the same, modern algebraic re-writing can, and here does, hide some difficulties. To better understand the difference of the points of view, let us note that the modern product $m \times n$ would be the number of squares in a grid-pattern with $m$ columns and $n$ rows (or more scholarly the number of elements of a Cartesian product of two sets respectively with $m$ and $n$ elements). The property of commutativity (i.e. $m \times n = n \times m$) is then evident. However, this is much less clear according to the definition of a product as successive additions of an integer to itself, for instance that the result of the addition of 3 to itself (3 + 3) is the same as twice the addition of 2 to itself (2 + 2 + 2). The asymmetry in the Greek definition of multiplication through additions is still clear in much later texts. For instance, what we consider as the products of two integers as 1×2, 2×3, 3×4, … are written by the 2$^{nd}$ century CE author Theon of Smyrna, respectively as: 'once 2' ('ἅπαξ μὲν γὰρ β´´'), twice 3 ('δὶς δὲ γ´´'), 'thrice 4' ('τρὶς δὲ δ´´'), and so on (Theon (1978), Arithmetic, XIII and also, X, XI, XII, …). In the same way, for a product, Jean Itard writes the integer as a capital letter, and the number of times it the product in lowercase: for instance $a$ times $b$ is written $Ba$ (Itard (1961), p. 104).

[6] For instance $4 = 1 \times 4$.

[7] Such precaution seems to have been forgotten in later times, cf. *infra*, note 12.

[8] Cf. Szabó (1978), p. 41, note 24 including Cherniss' remark.

[9] For instance, writing $3 \times 2 = 6$, the integer 6 is a 'plane integer', and its sides are 3 and 2.

[10] That may seem to follow, for arbitrary integer, from Theodorus' constructions. But, contrary to Theodorus, Theaetetus did not use the foot as unit; thus, though he did not specify it, we need to suppose that an

arbitrary (finite) line *u* had been set as the unit, so that any integer *n* is represented by a line equal to *n* times *u*, probably already such a common construction in the fifth century that it might go without saying.

[11] We refer to Brisson-Ofman (2017) for more details.

[12] Jean Itard translates also 'προμήκη' by 'rectangular' ('*rectangulaires*', Itard (1961), p. 216); nevertheless, at the beginning of his book, he considers that Euclid's 'plane number' is rather a 'rectangular number', thus including in it, contrary to Theaetetus, perfect squares (ib., p. 21). It is certainly not the only ambiguities with these terms, from the Antiquity to our days. In the mathematical passage, 'προμήκης' and 'ἑτερομήκης' are used indifferently (cf. *infra*, point 10)). In Aristotle's Pythagorean table of contraries, 'ἑτερομήκης' is opposed to 'square', thus has the same meaning as in the Theaetetus account, integers that can be written only as the product of two different integers (*Metaphysics*, I, 986a26). Knorr translates it either by transliteration as 'heteromecic' or, as 'oblong' (Knorr (1975), p. 144). Euclid does not use these terms in arithmetic (though in definition 22 of book I for instance, he considers 'ἑτερόμηκες' quadrilateral ('τετράπλευρος') defined as a (strict) rectangle i.e. that is not a square). However, they are found again in texts of late Antiquity. But their authors differentiate between 'προμήκη' and 'ἑτερομήκης', though they do not agree exactly on their meanings (cf. Heath (1908), vol. 2, p. 289). For instance, for Nicomachus, Theon of Smyrna as well as Iamblichus, 'ἑτερομήκης' means the kind of integer which has the power to be a product of two successive integers (i.e. of the form $n(n+1)$). However, while for Theon, 'προμήκη' indicates an integer product of two unequal terms (of the form $n(n+k)$, with *k* an integer, in agreement with Theaetetus definition, cf. Theon (1978), Arithmetic, XVII and XII), Nicomachus excludes the 'ἑτερομήκης' from it (thus, it means an integer of the form $n(n+k)$, with *k* greater or equal to 2). J. Dupuis translates 'προμήκη' and 'ἑτερομήκης' by transliteration ('*promèque*' and '*hétéromèque*', (1892), Arithmetic, XII, XVII), R. an D. Lawlor respectively as 'oblong' and 'unequilateral' (Theon (1978), ib.), J. Delattre Biencourt as 'elongated rectangle' and 'just-rectangle' ('*rectangle allongé*' and '*juste rectangle*', (2010), p. 140). Moreover, it should be noticed that the precautions used by Theaetetus to define both the 'square' and the 'rectangular' integers seem forgotten in Theon's time. *Stricto sensu*, any square integer is 'ἑτερομήκης', since any integer may be written as the product of the unit by itself (cf. *supra*, note 6)), so that, according to Theon's definition, it is 'ἑτερομήκης' (cf. *supra*, points 3) and 6); it cannot though be 'προμήκης' since no perfect square can be written as the product of two successive integers). Conversely, Theon rightly remarks that, according to his definitions, an integer may both be a product of two successive integers and of two integers differing by more than one, giving the example of 12 which is both thrice 4 (thus 'προμήκης') and twice 6 (thus 'ἑτερομήκης' and non-'προμήκης', i.e. a 'παραλληλόγραμμος' integer defined previously (Dupuis (1892), Arithmetic, XIV)). To top all this confusion, there is, to quote Itard ((1961), p. 21), the 'catastrophic' mistake in the proofs of propositions VIII.22 and 23 in Euclid's *Elements*, that is rooted in the forgetting that a square integer can be also rectangular (or 'προμήκη'). Since these propositions, are clear corollary of the results of book VII (cf. *infra*, **Appendix**, 3.a)), there is here a real possibility that this part of the *Elements* was modified in a late edition less careful than the mathematicians of the 5[th]-4[th] century BCE. Let us also remark lastly that the definition of such integers remains ambiguous still in modern time; for instance,

'rectangular number' means sometimes any non-prime integer (e.g. Higgins (2008), p. 9), sometimes a product of two consecutive integers (e.g. Ben-Menahem (2009), p. 161).

[13] Cf. *supra*, point 9) and 10) and also *supra*, note 12.

[14] Since two integers are (by their very definition) always commensurable, it is also equivalent to say that the side is commensurable either to the unit or to any given integer. Theodorus showed the first in his lesson.

[15] 'κύβος δὲ ὁ ἰσάκις ἴσος ἰσάκις ἢ [ὁ] ὑπὸ τριῶν ἴσων ἀριθμῶν περιεχόμενος.' (definition VII.19).

[16] In another terms, to any integer $n$ that is not a perfect cube, the parallelepiped figure of sides $(1,1,n)$ will be associated. Let us remark that as in the case of 'plane' integers, this geometrical figure is unique.

[17] Let us consider the integer 3. It is not a perfect cube, and the parallelepiped $(1,1,3)$ is associated to it (cf. previous note). Thus to this parallelepiped is associated a cube of the same volume (i.e. 3) so that, in modern writing the side of this cube is $\sqrt[3]{3}$. This magnitude $\sqrt[3]{3}$ is not commensurable 'as length' to the first one (defined in 1)a)v), *supra*) i.e. the integers, however it is commensurable 'as the cube it has the power to produce' since $\sqrt[3]{3} \times \sqrt[3]{3} \times \sqrt[3]{3} = 3$ is commensurable, as an integer, to any other integer (cf. *supra*, note 14).

[18] This is consistent with the former plane construction in the following sense: in any plane contained in the 3-dimensional space, the spatial 'length' is the same as the plane 'length' defined previously by Socrates-Theaetetus.

[19] Nevertheless, let us remark two pitfalls of their definition. The first is the same as for the plane case: in modern terms, the proof that the cubic root of any non-cube integer is incommensurable to all the integers is missing. The second is connected to the geometrical construction. Since in it is impossible to go further than the 3-dimensionnal space, there is no possible definition for roots of higher degree. In particular, it is not possible through this construction to define a square root of a square root i.e. a root of order 4 (for instance $\sqrt[4]{3}$), thus to obtain the geometric mean between, for instance, the unit and the square root of 3 (the line $x$ such that $1/x = x/\sqrt{3}$). Yet it is possible to construct such a line geometrically, using the same corollary of Pythagoras' construction that Theodorus used to draw his 'powers' in his lesson (cf. B-O, **Appendix I**).

[20] '[F]or number is a plurality of units' ('ὁ δ' ἀριθμὸς πλῆθος μονάδων', *Metaphysics*, 10, 1053a20); later he gives another one: 'number denotes a measured plurality and a plurality of measures' ('ὁ ἀριθμὸς ὅτι πλῆθος μεμετρημένον καὶ πλῆθος μέτρων', *ib.*, 14, 1088a5).

[21] See Dedekind (1961), in particular p. 6-12 for the modern comparison 'between numbers and points of a straight line', and p. 12-19 for the creation of 'numbers' in the modern sense.

[22] 'ἄπειροι τὸ πλῆθος αἱ δυνάμεις ἐφαίνοντο, πειραθῆναι συλλαβεῖν εἰς ἕν.' (147d7-8).

[23] For example, *Hippias* 287e-288a, *Meno* 72a-b. Since these dialogues belongs to Plato's early works, it is a testimony that he never changed his mind on this question.

[24] The problem is that according to many modern historians of mathematics, the unit was not considered as an integer by the ancient Greeks. Such a position is grounded in some Aristotelian texts (cf. *supra*, note 20), and in Euclid's definition of an integer as 'multitude composed of units'. Nevertheless, even in Euclid's *Elements*, the unit is often treated as an integer (e.g. proposition VII.15, proposition VIII.21 *vs.* definition VII.16), so that

Jean Itard can write that, in the treatise, 'the unit (…) is sometimes distinguished from number, sometimes identified to it' (Itard (1961), p. 71). In Plato's texts, including of course *Theaetetus*' mathematical passage, the unit is usually considered as an integer (e.g. *Great Hippias*, 302a; *Parmenides* 153a,e). This is not limited to early Antiquity, for example while Theon of Smyrna explicitly writes the unit is not an integer (VII, l. 29-30, adding by the way neither is the 'dyad'), he considers, and quotes authors that consider the unit belonging to the integers (III, l. 29-30; V, l. 10-11; VI, l. 21-22; VI, l. 1-2; 15-16; XIII, l. 1 and 4-5; XIX, l. 29 and l. 1;…; cf. also Plutarch, *Moralia*, 1003F). For hesitations by modern historians, see for instance Knorr (1975), p. 9 *vs*. p. 15, 24, 146-153.

[25] For the pitfalls of the boys' construction with respect to Euclid's construction, we refer to Brisson-Ofman (2017).

[26] For a similar analysis, see Knorr (1975), p. 84.

[27] Cf. *supra*, 2.1)14).

[28] As a matter of fact, according to the beginning of the passage, Theodorus did not prove such a result, since he stopped at 17. However, the boys say that they did not see any reason why he stopped there (cf. 147d6 and the translation in B-O). Hence, they naturally concluded that the result proved in the lesson for the first integers can be extended to any integer, as made clear in their later definition of the 'powers' (cf. previous note).

[29] In modern terms 'irrational numbers'.

[30] Since the square root of a perfect square is an integer and conversely any integer is the square root of this squared integer; in symbolic terms, for any integer *n*, we have: $\sqrt{n^2} = n$.

[31] It can be understood as a graphical representation to obtain general geometrical mean, as Aristotle writes in 'On the Soul': 'What is squaring? The construction of an equilateral rectangle equal to a given oblong rectangle. (…) squaring is the discovery of a line which is a mean proportional between the two unequal sides of the given rectangle discloses the ground of what is defined.' (*De Anima* II, 2, 413a16-19, J. Smith translation).

[32] Theaetetus switches from the singular to the plural

[33] In modern language, a power is a square roots of the form $\sqrt{n}$, with *n* an integer which is not a perfect square, so that its square is $(\sqrt{n})^2 = n$ is an integer, thus commensurable to any integer, hence the 'length'.

[34] In Brisson-Ofman (2017), we tried to show that Plato wanted to emphasize Theaetetus' concern about sensible things rather than their mutual connections. This is consistent with the replacement of the usual dative form (i.e. 'in power') as found in Euclid's and later mathematical texts by a name, 'power'. If our analysis is correct, the use of 'power' as a name instead of a relation, is a Theaetetus' idiosyncrasy and not a change in the mathematical language of ancient Greek mathematicians.

[35] Cf. also *supra*, note 23.

[36] Cf. §2), *supra*.

[37] Or in symbolic language, for *n*, *k* integers, we have the following equivalence:
$\sqrt{n} = k \Leftrightarrow n = k^2$ (called '**equivalence number 1**' in the **Appendix** 1.i)).

[38] Or in symbolic language:

$\sqrt{n}/k = p/q$ (*n, k, p* and *q* integers) ⇔ *n* is a perfect square (the '**equivalence number 2**' in the **Appendix 1.ii)**).

[39] Cf. *supra*, note 2.

[40] Or according to the Appendix, they think to '**equivalence number 1**' but speak of the '**equivalence number 2**' that requires a development of the theory of 'relatively prime integers' not available at the time of the lesson (see *supra*, Appendix 3.a), in particular note 89). The same gap between these two equivalences appears another time in the final generalization to cubic integers (148b2, cf. *supra*, point 2.1)15)).

[41] It is the reason why Pappus praises Euclid over Theaetetus (cf. *supra*, note 82).

[42] Pappus gives the example of the square roots of 2 and 8, which are commensurable to one another (for $\sqrt{8}/\sqrt{2} = 2$), but are outside of the classes considered by Theaetetus (Thomson (1930), §17, p. 82).

[43] For a detailed analysis concerning this problem, we refer to Brisson-Ofman (2017).

[44] See for instance Ryle (1966) and Owen (1965). Such claim has some consequences about the dating of the dialogues and also of Theaetetus' death explaining the success of its late dating among modern philosophers.

[45] G. Ryle considered the digression to be 'philosophically pointless' (Ryle (1966), p. 158). Mc Dowell considers that the digression 'is quite irrelevant to the dialogue' (McDowell (1973), p. 174). For Myles Burnyeat, in the digression 'Plato interrupts the argument and launches into rhetoric' (Burnyeat (1990), p. 34). According to Andrew Barker, the digression is to be explained merely by the desire to show what Protagoras' doctrine would lead to (Barker (1976), p. 462). Timothy Chapell opts for the following middle path: the 'Digression paints a picture of what it is like to live in accordance with two different accounts of knowledge, the Protagorean and the Platonist, that Plato is comparing (Chapell (2004), p. 127). Thus, the Digression shows us what is ethically at stake in the often abstruse debates found elsewhere in the *Theaetetus*. Its point is that we can't make a decision about what account of knowledge to accept without making all sorts of other decisions, not only about the technical, logical and metaphysical matters that are to the fore in the rest of the *Theaetetus*, but also about questions of deep ethical significance. So, for instance, it can hardly be an accident that, at 176c3, the difference between justice and injustice is said to be a difference between knowledge (*gnôsis*) and ignorance (*agnoia*)'. These positions depend on unstated presuppositions, as does this declaration by M. Burnyeat on the Forms: 'the question whether Forms are referred to at 174b and 175c is the question whether, *in those contexts*, the operative phrases are naturally taken to imply more metaphysics than (Aristotle) is actually expressed' (Burnyeat (1990), p. 38).

[46] See *supra*, I.2.16).

[47] Or companion, or disciple ('συγγιγνόμενος' 150d2, 'συνουσία', 150d3).

[48] The so-called 'recollection' of the servant-boy with regard to doubling a square of side 2-foot (82b-85b).

[49] Cf. Brisson (2007, 2008).

[50] Cornford (1935), p. 86, n. 1; *contra* Robinson (1950); response in Hackforth (1957).

[51] For detail, see Brisson-Ofman (2017).

[52] Translation by G.M.A Grube, revised by C.D.C. Reeve.

[53] 'So that I am not in any sense a wise man' ('εἰμὶ δὴ οὖν αὐτὸς μὲν οὐ πάνυ τι σοφός', 150d1; already claimed in 150c7-8).

[54] Terms relative to maieutics are rather frequent throughout the *Theaetetus*. The art of delivery: μαιευτική 149b6, c1, 150b6, 8, c7, c9, e5, 151c1, 160e3, 161e5, 184b1, 210b8; to conceive: κυεῖν 184b1, 210b4, 149c6; to give birth: τίκτειν 149b7, d2, 8, 150b8, 151a6; to suffer the pains of childbirth: ὠδίνειν 148e5, 149c6, d6, 151a6, 8, b8, 210b4, 150b1; 151b8

[55] E.g. *Theaetetus*, 161a7-b1; *Symposium* 176c3-e2.

[56] On the process of sense perception see Luc Brisson (1997).

[57] The advantages are interesting enough to be used even by contemporary mathematicians e.g. Nielsen (1993).

[58] See *infra*, **Appendix**.

[59] At *Timaeus* 64b, sensible data are sent to the intellect, which can make a judgment. See Brisson (1999; 2013).

[60] In Greek, the term '*dynamis*' means 'power' but also (human) 'faculty'.

[61] 'πάντα τὰ πρόσθεν ἐξαλείψας' (187b1).

[62] The term '*doxa*' ('δόξα') is given differently, sometimes in the translation of a same dialogue: 'opinion' by A. Diès (1926) and M. Narcy (1994), as opposed to 'judgment' by L. Robin (1950) for the French translations, judgment by J. McDowell (1973) and M. Burnyeat (1990) for the English. Often the translator chooses different words depending on the background, or at least provides a warning about the ambiguity of its meaning. Still another proposed translation is 'belief' (e.g. Rowe (2007), p. 229; also sometimes Chappell (2004)).

[63] *Phaedrus* (273b-c).

[64] *Phaedrus*, translated by A. Nehamas and P. Woodruff, Indianapolis, 1995.

[65] cf. *supra*, I.2.1)7).

[66] One of the most important problem is to get different sets of elements, as it is made clear in the mathematical passage: 6 is either 'after 1, 2, 3, 4, 5', or 'twice three', or 'three times 2', or '1 + 2 + 3' (204b9-c2). As a matter of fact, the correct proof of the boys' pseudo-definition is founded onto the notion of 'relatively prime' (in Greek 'πρῶτοι πρὸς ἀλλήλους', 'prime to one another', Euclid's *Elements*, definition VII.12) integers (cf. **Appendix**, §1.ii)), as opposed to (absolute) 'prime integer' ('πρῶτος ἀριθμός', *ib*, definition VII.11). For a general discussion about objects vs. ratios of objects in the *Theaetetus*, we refer to Brisson-Ofman (2017).

[67] According to the construction of 'rectangular integers' as the product of 1 by the integer itself, it would mean we know nothing both about the unit and the integer, but we know the integer.

[68] As a matter of fact, Socrates delivers Theaetetus not only of his present bad children (ideas) but also future ones, something a midwife cannot do.

[69] Cf. *supra*, §iii).

[70] Thus, contrary to the claim of many commentators based on a sentence Socrates uttered in the *Meno*, the correct definition of knowledge (for Plato) missing (correct) definition at the end of the *Theaetetus* cannot be 'true opinions/judgment fasted by causal reasoning'. As a matter of fact, Socrates never claimed that, 'true opinions/judgments (…) are of no great value until one makes them fast with causal reasoning' ('ἕως ἄν τις αὐτὰς δήσῃ αἰτίας λογισμῷ.', 98a4) may give the definition of knowledge.

[71] An assertion based on texts from late Antiquity supposed to account Plato's life and doctrine using second/third hand dubious compilations or even some forged texts. For instance, the supposed relation between Plato and Theodorus is mostly based on two sentences in passing by Diogenes Laertius. The first one appears in a list of twenty people named Theodorus which says that one of them was a geometer whose lectures were attending by Plato (Hicks (1925), II, 8, 103). The second in the life of Plato, where it is said that the philosopher had made a visit to Cyrene to Theodorus the mathematician (ib., III, 6, 103).

[72] Cf. *supra*, §I.4.

[73] Cf. previous note, and *supra*, 0.

[74] See *supra*, §0.3.iii).

[75] As a matter of fact, Gauss was definitively not the first to get this result, not even to state it. We find it already in a book of Jean Preslet printed in 1689: 'If an integer $d$ measures exactly the product $bc$ of two integers $b$ and $c$, and let $c$ and $d$ be prime to one another; the integer $d$ is a divisor of the other integer $b$' (Prestet (1689), I, book 6, corollary 3 of theorem 1, §21, p. 147; cf. Itard (1961), p. 113). While not stated in Euclid's *Elements*, it is a trivial consequence of propositions VII.20-21: starting with the same hypothesis, for some integer $k$, we get: $bc = kd$, thus $c/d = k/b$; since $c$ and $d$ are prime to one another, they are the least integers of ratio $k/b$ (proposition VII.21), thus (proposition VII.20) $d$ divides $b$ (*ib.*). QFD.

[76] Cf. *infra*, notes 80-82.

[77] Cf. *infra*, note 39.

[78] The definition of rationality/irrationality in Euclid's *Elements* is somewhat different from their usual meaning in Greek mathematics. To avoid any confusion, we will stick with the latter, since it has been the standard meaning from the Antiquity until today.

[79] We are quoting here the mathematician Barry Mazur who in (2008) pointed to the 'shocking unanimity' among commentators, qualifying it as a general 'strange delusion'. However, his interest was more about Euclid's proofs and the correction of a few mathematical defects, rather than in Plato's text.

[80] For instance Knorr (1983) p. 42: 'the general result incorporating these was known, and perhaps first enunciated by Theaetetus of Athens early in the fourth century. Plato gives a loose statement of it in the dialogue named after this mathematician: The mathematician'. He then quotes *Theaetetus* 148a-b to conclude: 'The sharp arithmetical cast of this formulation is tempered in the statement of the more general condition within the Euclidean theory' quoting proposition X.9 of the *Elements*, adding a little further 'Pappus rightly contrasts the restricted Platonic statement of this condition with the general Euclidean statement' (p. 44); cf. also Caveing (1998), p. 259-260; Burnyeat (1978), p. 509.

[81] Pappus of Alexandria was a mathematician of the 4th century CE. This commentary has been preserved only in an Arabic translation, which may add some difficulties to understand it (cf. for instance Knorr (1983), note 13, p. 63). We use here its English translation by Thompson (1930).

[82] 'This proposition is a Theaetetean discovery, and Plato has made mention of it in the *Theaetetus*; but there it occurs in a more restricted form, while here it is general.' ('Τὸ θεώρημα τοῦτο Θεαιτήτεόν ἐστιν εὕρημα, καὶ μέμνηται αὐτοῦ ὁ Πλάτων ἐν Θεαιτήτῳ, ἀλλ᾽ ἐκεῖ μὲν μερικώτερον ἔγκειται, ἐνταῦθα δὲ καθόλου.', scholium 62 to Euclid, *Elements* X ed. Heiberg, p. 450-452; cf., also Knorr (1975), p. 97, n. 6). We have to be cautious with this testimony, since chances are high it is just a remark quoted from Pappus' commentary. For a discussion on this proof, see Heath (1908), vol. 3, p. 29-31.

[83] The proof uses some results not given in Euclid's treatise, e.g. the equivalence that two ratios and the ratio of their squares are the same. Depending on the point of view concerning the general theory of proportions established in book V, this is either a triviality (cf. next note) or something almost impossible to prove. The difficulty (if any) would come from the combination of ratios of magnitudes and ratios of integers. For the analysis of the proof, see for instance, Heath (1908), vol. 3, p. 30-31 and Caveing (1998), p. 236-245. In Theaetetus' account, this equivalence is entailed by the definition itself.

[84] In the case considered by Theaetetus, this equivalence is given by proposition VIII.11 of Euclid's *Elements*, but was known much earlier, for it is needed to show the incommensurability of the diagonal by the theory of odd and even, as presented it in Aristotle's *Analytics* (cf. Ofman (2010), p. 114; Ofman (2014), **Appendix II**). Pappus emphasizes that Proposition X.9 is more general since it deals not only with integers but general magnitudes (cf. *supra*, note 82).

[85] The algebraic writing is misleading since the difficulty in the long proof of proposition X.9 is not the squaring of ratios. As a matter of fact, in Euclid's *Elements*, this is considered as self-evident and not proved within the treatise. Indeed, there is not even a definition in an arithmetical context of a 'squared' ratio, while such a definition is given for general magnitudes (definition V.9; for magnitudes, Euclid uses the term '*diplasiona logon*' i.e. 'double ratio'). These are more testimonies that this part of the *Elements* follows some old well-established theory. This is the common view among historians as different as e.g. Zeuthen (1910) Heath (1908), Becker (1936), van der Waerden (1963) and Knorr (1975) who considers this question in detail (p. 239-244). The difficulty is rather lying in the inhomogeneity of the terms in ratio, two magnitudes on the one hand, two integers in the other hand.

[86] For $\sqrt{n}/m = p/q$ is equivalent to: $\sqrt{n} = mp/q$. For the questions about commensurability/incommensurability and rationality/irrationality, we refer to B-O, Appendix III.

[87] In modern sense, a rational 'number', or briefly a 'rational', is a real 'number' which can be written as the quotient *p/q* of two integers. Putting in '**equivalence number 3**': $r = c/d = \sqrt{a}/\sqrt{b}$, we use the identity: $\sqrt{a}/\sqrt{b} = \sqrt{a/b}$ (cf. *supra*, note 83).

[88] Definition VII.13: 'Numbers prime to one another are those which are measured by an unit alone as a common measure.' i.e. they have no common divisor (except 1). For instance, 3 and 10 are relatively prime while 4 and 10 are not relatively prime since 2 divides both.

[89] We refer to B-O, more precisely **Appendix** 2, where we studied the proof of the (spurious) proposition X.117 in the *Elements*.

[90] 'If $c$ and $d$ are the least integers such that $a/b = c/d$ then $c$ (respectively $d$) divides $a$ (respectively $b$).' (proposition VII.20); these least integers are the ones which are prime to one another. (Propositions VII.21-22).

[91] Since it is evident that one of the primary purposes of propositions VII.20-27 is this result (cf. also, *supra*, note 12, this reconstruction is hardly the first along these lines (e.g. van der Waerden (1963), p. 168; Knorr (1975), p. 229).